\documentclass[11pt]{article}

\usepackage[latin9]{inputenc}
\usepackage[T1]{fontenc}
\usepackage{pifont}
\usepackage{makeidx}
\usepackage{graphics}
\usepackage{amsmath,latexsym}
\usepackage[psamsfonts]{amssymb}

\newtheorem{theorem}{Theorem}[section]
\newtheorem{prop}[theorem]{Proposition}

\newtheorem{definition}[theorem]{Definition}
\newtheorem{question}[theorem]{Question}

\newtheorem{examples}[theorem]{Examples}

\newcommand\thmm[2]{\begin{theorem}[#2] #1. \end{theorem}}

\def\dfn#1{\textbf{#1}}  

\newcommand{\action}[1]{ {\curvearrowright}^{\! #1}}
\def\acting{\curvearrowright}

\def\RR{\mathcal{R}}
\def\SS{\mathcal{S}}
\def\SL{\mathrm{SL}}
\def\VV{\texttt{V}}
\def\Nmath{\mathbb{N}}
\def\Smath{\mathbb{S}}
\def\Cmath{\mathbb{C}}
\def\Rmath{\mathbb{R}}

\def\Zmath{\mathbb{Z}}
\def\Tmath{\mathbb{T}}
\def\Leb{\mathrm{Leb}}
\def\Aut{\mathrm{Aut}}
\def\Out{\mathrm{Out}}
\def\FF{\mathbf{F}}

\def\acting{ \curvearrowright }

\def \hyp{\mathrm{hyp}}
\def\Conj{\mathop{\overset{\textrm{\scriptsize{Conj}}}{\sim}}}
\def\OrbEqu{\mathop{\overset{\textrm{\scriptsize{OE}}}{\sim}}}
\def\vNEqu{\mathop{\overset{\textrm{\scriptsize{vN}}}{\sim}}}
\def\ME{\mathop{\overset{\textrm{\scriptsize{ME}}}{\sim}}}
\def\MFI{$\mathcal{MFI}$ }

\def\cost{\mathrm{Cost}}

\def\geomdim{\textrm{geo-dim}}
\def\approxdim{\textrm{approx-dim}}

\begin{document}

\title{Orbit Equivalence and Measured Group Theory}
\author{Damien Gaboriau\thanks{CNRS / ANR AGORA N°ANR-09-BLAN-0059}}
\date{May, 2010}

\maketitle

\begin{abstract} We give a survey of various recent developments in orbit equivalence and measured group theory. This subject aims at studying infinite countable groups through their measure preserving actions.
\end{abstract}

\bigskip
\noindent
{\small \textbf{2000 Mathematical Subject Classification}: Primary 37A20; Secondary 46L10; 28D15.\\
\noindent
\textbf{Key words}:  Orbit equivalence, Measured group theory, von Neumann algebras.}

\section{Introduction}
Orbit equivalence and measure equivalence theories deal with countable groups $\Gamma$ acting on standard measure spaces 
 and with the associated orbit partitions of the spaces.
This is very much connected from its birth with operator algebras \cite{MvN36}; many of the recent progresses in both areas were made conjointly (see \cite{Popa=ICM-2007,Vaes=Bourbaki=2007,Vaes=ICM=2010}). 
It turns out to be also connected with geometric group theory (see section~\ref{sect:Measure equivalence} and \cite{2009=Furman-survey}), descriptive set theory (see \cite{JKL02,Kechris-Miller}), percolation on graphs (see \cite{Lyons-Peres-book})... with fruitful cross-pollination. 

There are many examples of mathematical domains where the orbit equivalence or measured approach helps solving delicate questions involving countable groups $\Gamma$. For instance, in connection with group $\ell^2$-Betti numbers $\beta^{_{(2)}}_n(\Gamma)$, this was useful
to attack:\\
-- various vanishing results in \cite{Gab02,Sauer-Thom=spectral-seq=L2=2007};\\
-- the study of harmonic Dirichlet functions on percolation subgraphs \cite{Gab05};\\
-- the comparison between the uniform isoperimetric constant and $\beta^{_{(2)}}_1(\Gamma)$ \cite{Lyons-Pichot-Vassout-08};\\
-- problems of topological nature, related to the work of Gromov about the minimal volume \cite{Sauer=Amenable-cover+volume+L2=2009}.

In geometric group theory, the quasi-isometry invariance of various cohomological properties for amenable groups \cite{Shalom04} was obtained that way.
Gaboriau-Lyons' measurable solution to von Neumann's problem \cite{Gaboriau-Lyons-2009} happens to be a way to extend results about groups containing a copy of the free group $\FF_2$ to every non-amenable group (see section~\ref{sect:How many actions}). This was used by \cite{Epstein-07} and in Dixmier's unitarizability problem \cite{Monod-Epstein=2009,Monod-Ozawa=Dixmier=2010}. 

The purpose of this survey is to describe some foundations of the theory and some of its most recent developments. There are many aspects upon which we shall inevitably not touch here, and many results are just alluded to, with as far as possible the relevant bibliography.

There are several excellent books and surveys with various focuses on orbit equivalence to which the reader is referred for further information, for instance \cite{Kechris-Miller,Gab05a,Shalom=survey=2005,Popa=ICM-2007,2009=Furman-survey,Kechris=Global-aspects=2010}.

\section{Setting and examples}

 The measure spaces $X$ will always be assumed to be standard Borel spaces and unless specified otherwise, the measure $\mu$ will be a non-atomic probability measure. 
Measurably, $(X,\mu)$ is isomorphic to the interval $([0,1], \Leb )$ equipped with the Lebesgue measure.
Moreover, the actions $\Gamma\action{\alpha}(X,\mu)$ we shall consider will be by Borel automorphisms and probability measure preserving (p.m.p.), i.e. $\forall \gamma\in \Gamma$, $A\subset X$: $\mu(\gamma.A)=\mu(A)$ (We only considers Borel sets).
Shortly, $\alpha$ is a \dfn{p.m.p. action} of $\Gamma$.
In this measured context, null sets are neglected. Equality for instance is always understood almost everywhere. The action $\alpha$ is  \dfn{(essentially) free} if $\mu\{x:\ \gamma.x=x\}>0 \Rightarrow \gamma=id$. 
The action is \dfn{ergodic} if the dynamics is indecomposable, i.e. whenever $X$ admits a partition $X=A\cup {^c\!{A}}$
into invariant Borel subsets, then one of them is trivial, i.e. $\mu(A) \mu({^c\!{A}})=0$.

We now present a series of basic examples which shall already exhibit a rich variety of phenomena.
\begin{examples}\label{ex: action on the circle}
The action of $\Zmath^{n}$ on the circle $\Smath^{1}$ by rationally independent rotations.
\end{examples}

\begin{examples}\label{ex: SLn acting on the n-torus}
 The standard action $\SL(n,\Zmath)\acting\Tmath^n$ on the $n$-torus $\Rmath^n/\Zmath^n$ with the Lebesgue measure. The behavior is drastically different for $n\geq 3$ and for $n=2$. The higher dimensional case was central in the super-rigidity results of Zimmer \cite{Zimmer-Erg-th-ss-gp-Book-84} and Furman \cite{Fur99a, Fur99b} (see section~\ref{sec:some rigidity results}).
 The $2$-dimensional case $\SL(2,\Zmath)\acting\Tmath^2$ played a  particularly important role in the recent developments of the theory, mainly because of its relation with the semi-direct product $\SL(2,\Zmath)\ltimes \Zmath^2$, in which $\Zmath^2$ is known to have the so called relative property (T) (see section~\ref{sec:Rel T}), 
 while $\SL(2,\Zmath)$ is a virtually free group (it has a finite index free subgroup).
\end{examples}

\begin{examples}
Volume-preserving group actions on finite volume manifolds.
\end{examples}

\begin{examples} Given two lattices $\Gamma, \Lambda$ in a Lie group $H$ (or more generally a locally compact second countable group) the actions by left (resp. right by the inverse) multiplication on $H$ induce actions on the finite measure standard spaces $\Gamma \acting H/\Lambda$ and $\Lambda \acting \Gamma\backslash H$.
\end{examples}

\begin{examples} A compact group $K$, its Haar measure $\mu$  and the action of a countable subgroup $\Gamma$ by left multiplication
on $K$.
\end{examples}

\begin{examples} \label{ex: Bernoulli shifts}
Let $(X_0,\mu_0)$ be a standard probability measure space, possibly with atoms\footnote{for instance $X_0=\{0,1\}$ and $\mu_0(\{0\})=1-p, \mu_0(\{1\})=p$ for some $p\in (0,1)$. The only degenerate situation one wishes to avoid is $X_0$ consisting of one single atom.}.
 The \dfn{standard Bernoulli shift action} \index{Bernoulli @! shift} \index{action @! Bernoulli shift} of $\Gamma$ is the action on the space $X^{\Gamma}$ of sequences $(x_{\gamma})_{\gamma\in \Gamma}$ by shifting the indices $g.(x_{\gamma})_{\gamma\in \Gamma}=(x_{g^{-1}\gamma})_{\gamma\in \Gamma}$, together with the $\Gamma$-invariant product probability measure $\otimes_\Gamma \mu_0$.
In particular, every countable group admits at least one p.m.p. action.
The action is {free} (and ergodic) iff $\Gamma$ is infinite. 

More generally, consider  some action $\Gamma\acting \VV$ of $\Gamma$ on some countable set $\VV$.
 The \dfn{generalized Bernoulli shift action} \index{Bernoulli @! generalized '' shift} \index{action @! generalized Bernoulli shift} of $\Gamma$ is the action on the space $X^{\VV}$ of sequences $(x_{v})_{v\in \VV}$ by shifting the indices $g.(x_{v})_{v\in \VV}=(x_{g^{-1}.v})_{v\in \VV}$, with the invariant product probability measure. 
\end{examples}

\begin{examples} \dfn{Profinite actions}. \index{profinite action} \index{action @! profinite}
Consider an action $\Gamma\acting(\mathtt{T},v_0)$ of $\Gamma$ on a locally finite rooted tree.
The action preserves the equiprobability on the levels, and the induced limit probability measure on the set of ends of the tree is $\Gamma$-invariant.
For instance, if $\Gamma$ is residually finite, as witnessed by a chain of finite index subgroups $\Gamma=\Gamma_0>\Gamma_1>\Gamma_2>\cdots\Gamma_i>\cdots$ with trivial intersection, 
such a rooted tree $(\mathtt{T}, (v_0=\Gamma/\Gamma_0))$ is naturally built with vertex set (of level $i$) the cosets $\Gamma/\Gamma_i$ and edges given by the reduction maps $\Gamma/\Gamma_{i+1}\to \Gamma/\Gamma_{i}$. The action is ergodic iff it is transitive on the levels.
\end{examples}

A first connection with functional analysis is made through the following.
The \dfn{Koopman representation} of a p.m.p. action $\Gamma\action{\alpha}(X,\mu)$ is the representation $\kappa_{\alpha}$ of $\Gamma$ on $L_0^2(X,\mu)$ given by\footnote{The constants are fixed vectors for the representation on $L^2(X,\mu)$. Its orthocomplement $L_0^2(X,\mu)=L^2(X,\mu){\ominus} \Cmath 1$ consists in $\{\xi\in L^2(X,\mu): \int_X \xi(x) d\mu(x)=0\}$.} $\kappa_{\alpha}(\gamma)(\xi)(x)=\xi(\alpha(\gamma^{-1})(x))$ \cite{Koopman-1931}.

A lot of dynamical properties of the action are read from this unitary representation and its spectral properties.
For instance, the action is ergodic if and only if its Koopman representation has no $\Gamma$-invariant unit vector.
In examples~\ref{ex: action on the circle} and~\ref{ex: SLn acting on the n-torus} or~\ref{ex: Bernoulli shifts}, various properties are deduced from the fact that the Koopman representation admits a Hilbert basis which is either made of eigenvectors or permuted by $\Gamma$ (see for instance \cite{Sch80}, \cite{Kechris-Tsankov-2008}). 
The classical ergodic theory considers such actions up to \dfn{conjugacy} (\textbf{notation:} $\Gamma_1\action{\alpha_1}X_1 \Conj \Gamma_2\action{\alpha_2}X_2$).

We now introduce a weaker notion of equivalence and turn from classical ergodic theory to orbit equivalence theory. Here $\Gamma.x$ denotes the orbit of $x$ under the $\Gamma$-action.
\begin{definition}[Orbit equivalence]\index{Orbit equivalence @!}
Two actions $\Gamma_i\action{\alpha_i}(X_i,\mu_i)$ (for $i=1,2$)
are \dfn{orbit equivalent (OE)} (\textbf{notation:} {$\Gamma_1\action{\alpha_1}X_1 \OrbEqu \Gamma_2\action{\alpha_2}X_2$})
if there is a measured space isomorphism\footnote{An isomorphism of measure spaces is defined almost everywhere and respects the measures: $f_*\mu_1=\mu_2$. \label{footnote=measure-space-isomoprhims=measure}} $f:X_1\to X_2$ that sends orbits to orbits:
\begin{center}
for a.a. $x\in X_1$: $f(\Gamma_1 .x)=\Gamma_2.f(x)$.
\end{center}
\end{definition}
In particular, the groups are no longer assumed to be isomorphic.
When studying actions up to orbit equivalence, what one is really interested in, is the partition of the space into orbits or equivalently the \dfn{orbit equivalence relation}:
\begin{equation}
\RR_{\alpha}:=\{(x,y)\in X: \exists \gamma\in \Gamma \textrm{\  s.t.\  } \alpha(\gamma)(x)=y\}.
\end{equation}

This equivalence relation satisfies the following three properties: (1) its classes are (at most) countable, (2) as a subset of $X\times X$, it is measurable, (3) it \dfn{preserves the measure} $\mu$: this means that every measurable automorphism $\phi: X\to X$ that is \dfn{inner} ($x$ and $\phi(x)$ belong to the same class for a.a. $x\in X$) has to preserve $\mu$.

Axiomatically \cite{FM77a}, the object of study is an equivalence relation $\RR$ on $(X,\mu)$ satisfying the above three conditions: we simply call it a \dfn{p.m.p. equivalence relation}.
Two p.m.p. equivalence relations $\RR_1, \RR_2$ will be \dfn{orbit equivalent}  if there is a measured space isomorphism $f:X_1\to X_2$ sending classes to classes.

This abstraction is necessary when one wants to consider, for instance, the \dfn{restriction} $\RR\vert A$ of $\RR$ to some non-null Borel subset $A\subset X$: the standard Borel space $A$ is equipped with the normalized probability measure $\mu_A(C)=\mu(C)/\mu(A)$ and 
$(x,y)\in \RR\vert A \Leftrightarrow x,y\in A \textrm{ and } (x,y)\in \RR$.

\medskip
In fact, this more general context allows for much more algebraic flexibility since the lattice of subrelations of $\RR_\alpha$ for some $\Gamma$-action $\alpha$ is much richer than that of subgroups of $\Gamma$ (see von Neumann's problem in section~\ref{sect:How many actions}). Also, $\RR_\alpha$ is easier to decompose as a ``free product or a direct product'' than $\Gamma$ itself (see section~\ref{sect: Dimensions} and \cite{Alvarez-Gaboriau-08}).

\medskip
By an \dfn{increasing approximation} $\RR_n\nearrow \RR$ of a p.m.p. equivalence relation $\RR$ we mean an increasing sequence of standard (p.m.p.) equivalence subrelations with $\cup_n\RR_n=\RR$.

An important notion is that of hyperfiniteness: a p.m.p. equivalence relation $\RR$ is \dfn{hyperfinite} if it admits an increasing approximation by \dfn{finite} equivalence subrelations $\RR_n$ (i.e. the classes of the $\RR_n$ are finite).
Obviously all the actions of locally finite groups (i.e. groups all of whose finitely generated subgroups are finite) generate orbit equivalence relations in this class; for instance such groups as $\Gamma=\oplus_{\Nmath} \Lambda_n$, where the $\Lambda_n$ are finite. This is also the case for all $\Zmath$-actions. 
Dye's theorem is among the fundamental theorems in orbit equivalence theory:
\thmm{All the ergodic hyperfinite p.m.p. equivalence relations are mutually orbit equivalent}{\cite{Dye59}} 
A series of results due in particular to Dye, Connes, Krieger, Vershik... leads to Ornstein-Weiss' theorem (see \cite{CFW81} for a more general version):
\thmm{\label{th: Ornstein-weiss} If $\Gamma$ is amenable then all its p.m.p. actions are hyperfinite}{\cite{OW80}}  In particular, when ergodic, these actions are indistinguishable from the orbit equivalence point of view! All the usual ergodic theoretic invariants are lost. This common object will be denoted $\RR_{\hyp}$.
On the other hand, if $\Gamma$ admits a free p.m.p. hyperfinite action, then $\Gamma$ has to be amenable, thus showing the border of this huge singular area that produces essentially a single object.
The non-amenable world is much more complicated and richer.

\section{The full group}
The \dfn{full group} of $\RR$ denoted by $[\RR]$ is defined as the group of p.m.p. automorphisms of $(X,\mu)$  whose graph is contained in $\RR$: $$[\RR]:=\{T\in \Aut(X,\mu): (x,T(x))\in \RR \textrm{\ for a.a.\ } x\in X\}.$$
It was introduced and studied by Dye \cite{Dye59}, and it is clearly an OE-invariant. But conversely, its algebraic structure is rich enough to remember the equivalence relation:
\thmm{(Dye's reconstruction theorem) Two ergodic p.m.p. equivalence relations $\RR_1$ and $\RR_2$ are OE iff their full groups are algebraically isomorphic; moreover the isomorphism is then implemented by an orbit equivalence}{\cite{Dye63}}
The full group has very nice properties. The topology given by the bi-invariant metric $d(T,S)=\mu\{x: T(x)\not= S(x)\}$ is Polish. In general, it is not locally compact and, in fact, homeomorphic with the separable Hilbert space $\ell^2$ \cite{Kittrell-Tsankov}.
\thmm{The full group is a simple group iff $\RR$ is ergodic}{\cite{Bezuglyi-Golodets=1980,Kechris=Global-aspects=2010}} 
And it satisfies this very remarkable, automatic continuity:
\thmm{If $\RR$ is ergodic, then every group homomorphism $f:[\RR]\to G$ with values in a separable topological group is automatically continuous}{Kittrell-Tsankov \cite{Kittrell-Tsankov}}

Hyperfiniteness translates into an abstract topological group property: 
\thmm{Assuming $\RR$ ergodic, 
$\RR$ is hyperfinite iff $[\RR]$ is extremely amenable}{Giordano-Pestov \cite{Giordano-Pestov=2007}}
Recall that a topological group $G$ is \dfn{extremely amenable} if every continuous action of $G$ on a (Hausdorff) compact space has a fixed point.
Together with Kittrell-Tsankov's result, this gives that every action of $[\RR_{\hyp}]$ by homeomorphisms on a compact metrizable space has a fixed point.

Closely related to the full group, the \dfn{automorphism group} $\Aut(\RR):=\{T\in \Aut(X,\mu): (x,y)\in \RR \Rightarrow (T(x),T(y))\in \RR\textrm{\ for a.a.\ } x\in X\}\triangleright [\RR]$ and the \dfn{outer automorphism group} (the quotient) $\Out(\RR)=\Aut(\RR)/[\RR]$ have attracted much attention for several years; see for instance \cite{1988=Gefter-Golodets=fund-gp,Gef93,Gef96,Fur05,IPP05,Popa-1-cohomology-2006,Kechris=Global-aspects=2010,Kida=Outer-2008,Popa-Vaes-2008,PV08,Gab08} and references therein and section~\ref{sec:Rel T}.

\section{Associated von Neumann algebra}
\label{sect: associated vN alg}
In fact, the original interest for orbit equivalence came from its connection with von Neumann algebras. Murray and von Neumann \cite{MvN36} considered p.m.p. group actions $\Gamma\action{\alpha} (X,\mu)$ as a machine to produce finite von Neumann algebras $M_{\alpha}$, via their group-measure-space construction. And Singer \cite{Singer-1955} was the first to explicitly notice that $M_{\alpha}$ only depends on the OE class of the action. Feldman-Moore \cite{FM77b} extended the group-measure-space construction to the context of p.m.p. equivalence relations.

A p.m.p. equivalence relation $\RR$ on $(X,\mu)$, considered as a Borel subspace of $X\times X$ is naturally equipped with a ({\textit{a priori}} infinite) \dfn{measure $\nu$}. It is defined as follows:
for every Borel subset $C\subset \RR$, 
\begin{equation}\label{eq: def nu}
\nu(C)=\int_X \vert \pi_l^{-1}(x)\cap C\vert d\mu(x),
\end{equation}
 where $\pi_l:\RR\to X$ is the projection onto the first coordinate, $\pi_l^{-1}(x)$ is the fiber above $x\in X$, and $\vert \pi_l^{-1}(x)\cap C\vert$ is the (at most countable) cardinal of its intersection with $C$.
A similar definition could be made with the projection $\pi_r$ on the second coordinate instead, but the fact that $\RR$ is p.m.p. ensures that these two definitions would coincide. 

The (generalized) \dfn{group-measure-space von Neumann algebra $L(\RR)$ associated with $\RR$} is generated by two families of operators of the separable Hilbert space $L^2(\RR,\nu)$:
$\{L_g: g\in [\RR]\}$ and $\{L_f:f\in L^{\infty}(X,\mu)\}$,  where $L_g\xi(x,y)=\xi(g^{-1}x,y)$ and $L_f\xi(x,y)=f(x)\xi(x,y)$ for every $\xi\in L^2(\RR,\nu)$.
It contains $\{L_f:f\in L^{\infty}(X,\mu)\}\simeq L^{\infty}(X,\mu)$ as a \dfn{Cartan subalgebra} (i.e. a maximal abelian subalgebra whose normalizer generates $L(\RR)$).
With this definition, $L(\RR)$ is clearly an OE-invariant.
\begin{definition}[von Neumann equivalence or W$^*$-equivalence]
Two p.m.p. equivalence relations  $\RR_i$ on $(X_i,\mu_i)$ (for $i=1,2$) are \dfn{von Neumann equivalent} or \dfn{W$^*$-equivalent}  if $L(\RR_1)\simeq L(\RR_2)$
(\textbf{notation:} $\RR_1 \vNEqu \RR_2$).
\end{definition}
There exist non-OE equivalence relations producing isomorphic $L(\RR)$ (\cite{CJ82}, \cite{ozawa-Popa-II=2008-b}).
Indeed, the additionnal data needed to recover $\RR$ is the embedding $L^{\infty}(X,\mu)\subset L(\RR)$ of the Cartan subalgebra  inside $L(\RR)$ (up to isomorphisms) \cite{Singer-1955,FM77b}.

\section{Strong ergodicity}\label{sect:strong ergodicity}

Recall that a standard  p.m.p. equivalence relation $\RR$ is \dfn{ergodic} if every $\RR$-invariant\footnote{for each $g$ in  the full group $[\RR]$, $\mu(A\Delta g A)=0$.} Borel set $A\subset X$  satisfies $\mu(A)(\mu(A)-1)=0$.
The notion of strong ergodicity was introduced by Schmidt as an OE-invariant.
\begin{definition}[\cite{Sch80}]\label{def:strong ergodicity}
An ergodic p.m.p. countable standard equivalence relation $\RR$ is \dfn{strongly ergodic} \index{strongly ergodic}\index{ergodic! strongly} if every almost invariant sequence\footnote{i.e. for each $g$ in $[\RR]$, $\lim_{n\to \infty} \mu(A_n\Delta g A_n)=0$.} of Borel subsets $A_n\subset X$ is trivial, i.e. satisfies $\lim_{n\to \infty}\mu(A_n)(1-\mu(A_n))=0$.
\end{definition}
There are several equivalent definitions of strong ergodicity, see for instance \cite{Jones-Schmidt-1987}.
We give yet another one below through approximations.
\begin{prop}\label{prop: def strong erg by approx}
An ergodic equivalence relation $\RR$ is \dfn{strongly ergodic} if and only if every increasing approximation $\RR_n\nearrow \RR$
admits an ergodic restriction $\RR_n\vert U$ to some  non-negligeable Borel set $U$, for big enough $n$.
\end{prop}
In other words, for big enough $n$ the ergodic decomposition of $\RR_n$ admits an atom.
It is easy to see that whenever a p.m.p. action $\Gamma\acting (X,\mu)$ is non-strongly ergodic, 
its Koopman representation $\kappa_0$ almost has invariant vectors.
 The converse does not hold in general \cite{Sch81}, \cite{Hjorth-Kechris=Rigid+products-05}. However, Chifan-Ioana \cite{Chifan-Ioana=Bernoulli-strong} extending an argument of Abert-Nikolov \cite{Abert-Nikolov-07} proved that this is indeed the case when the commutant of $\Gamma\acting (X,\mu)$ in $\mathrm{Aut}(X,\mu)$ acts ergodically on $(X,\mu)$. Standard Bernoulli shifts are strongly ergodic iff the group is non-amenable. In particular every non-amenable group admits at least one strongly ergodic action.

Kechris-Tsankov \cite{Kechris-Tsankov-2008} characterized the generalized Bernoul\-li shifts $\Gamma\acting (X_0,\mu_0)^{\VV}$ that are strongly ergodic as those for which the action $\Gamma\acting {\VV}$ is \dfn{non-amenable} (i.e. the representation on $\ell^2(\VV)$ does not admit any sequence of almost invariant vectors). 

The consideration of the Koopman representation $\kappa_0$ ensures that for (infinite) groups with Kazhdan property (T) \textbf{every} ergodic p.m.p. action is strongly ergodic. And Connes-Weiss  (by using Gaussian random variables) showed that this is a criterion for property (T) \cite{Connes-Weiss=1980}.

\medskip
A graphing $\Phi$ (see section~\ref{sect:cost}) on $X$ naturally defines a ``metric'' $d_\Phi$ on $X$: the simplicial distance associated with the graph structure in the classes of $\RR_\Phi$ and $d_\Phi=\infty$ between two points in different classes. 
This is a typical instance of what Gromov calls a \dfn{\textit{mm}-space} \cite{Gromov-00=spaces-questions}, i.e. a probability measure space $(X, \mu)$ together with a Borel function $d:X\times X\to \Rmath^{+}\cup\{\infty\}$ satisfying the standard metric axioms except that one allows $d(x,x')=\infty$.
A \textit{mm}-space $(X,\mu,d)$ is \textbf{concentrated} if
 $\forall \delta>0$, there is $\infty>r_\delta>0$ such that 
$\mu(A),\mu(B)\geq \delta \Rightarrow d(A,B)\leq r_\delta$.
For instance, if $\Phi=\langle \varphi_1:X\to X\rangle$ is given by a single p.m.p. ergodic isomorphism, $(X,\mu,d_{\Phi})$ is never concentrated.
Gromov observed for finitely generated groups that every p.m.p. ergodic action of $\Gamma$ has (respectively, never has) the concentration property if $\Gamma$ has Kazhdan's property (T) (respectively, if $\Gamma$ is amenable). 
Pichot made the connection with strong ergodicity:
\begin{theorem}[{\cite{Pichot-07=erg-forte}}]
Let $\Phi=(\varphi_i)_{i=1, \cdots, p}$ be a graphing made of finitely many partial isomorphisms. The space $(X,\mu,d_\Phi)$ is concentrated iff $\RR_{\Phi}$ is strongly ergodic.
\end{theorem}
See also \cite{Pichot=th-spectrale-rel-equ=2007} for a characterization of strong ergodicity (as well as of property (T) or amenability) in terms of the spectrum of diffusion operators associated with random walks on the equivalence relation $\RR$.

\medskip

For the standard $\SL(2,\Zmath)$ action on the $2$-torus $\Rmath^2/\Zmath^2$, every non-amenable subgroup $\Lambda<\SL(2,\Zmath)$ acts ergodically, and even strongly ergodically. Similarly for the generalized Bernoulli shift $\Gamma\acting (X_0,\mu_0)^{\VV}$, where the stabilizers of the action $\Gamma\acting \VV$ are amenable. Inspired by \cite{Chifan-Ioana=Bernoulli-strong}, define more generally:

\begin{definition}[Solid ergodicity]
A p.m.p. standard equivalence relation $\RR$ is called \dfn{solidly ergodic} if for every (standard) subrelation $\SS$ 
there exists a measurable partition $\{X_i\}_{i\geq 0}$ of $X$ in $\SS$-invariant subsets such that: \\
(a) the restriction $\SS\vert X_0$ is hyperfinite\\
(b) the restrictions $\SS\vert X_i$ are strongly ergodic for every $i>0$.
\end{definition}
In particular, an ergodic subrelation of a solidly ergodic relation is either hyperfinite or strongly ergodic.
By Zimmer \cite[Prop. 9.3.2]{Zimmer-Erg-th-ss-gp-Book-84}, every ergodic p.m.p. standard equivalence relation $\RR$ contains an ergodic hyperfinite subrelation $\SS$ which, being non strongly ergodic, contains an aperiodic subrelation with diffuse ergodic decomposition. Thus the $X_0$ part cannot be avoided, even for aperiodic subrelations.

One gets an equivalent definition if one replaces ``strongly ergodic'' by ``ergodic'' (see \cite[Prop. 6]{Chifan-Ioana=Bernoulli-strong} for more equivalent definitions).
It may seem quite unlikely that such relations really exist. 
However, Chifan-Ioana \cite{Chifan-Ioana=Bernoulli-strong}  observed that the notion of solidity and its relative versions introduced by Ozawa \cite{Ozawa=solidity=2004} (by playing between C$^*$- and von Neumann algebras) imply solid ergodicity (hence the name).
Moreover, they established a general solidity result for Bernoulli shifts.
\begin{theorem}The following actions are solidly ergodic:\\
-- The standard  action $\SL(2,\Zmath)\acting\Rmath^2/\Zmath^2$ \cite{Ozawa=2009-SL2xZ2-solid}.
\\
-- The generalized Bernoulli action $\Gamma\acting (X_0,\mu_0)^{\VV}$, when the $\Gamma$-action $\Gamma\acting {\VV}$ has amenable 
stabilizers \cite{Chifan-Ioana=Bernoulli-strong}.
\end{theorem}
When the group $\Gamma$ is exact\footnote{Recall that a discrete group $\Gamma$ is \dfn{exact} iff it acts amenably on some compact topological space.}, the above statement for the standard Bernoulli shifts also follows from \cite[Th. 4.7]{Ozawa-2006=Kurosh-type}.

A positive answer to the following percolation-theoretic question would give another proof of solid ergodicity for the standard Bernoulli shifts:
\begin{question}
Let $\Gamma$ be a countable group with a finite generating set $S$. 
Let $\pi:(X_0,\mu_0)^\Gamma \to [0,1]$ be any measure preserving map (i.e. $\pi_*(\otimes_{\Gamma} \mu_0)=\Leb$) and $\Phi_{\pi}$ be the ``fiber-graphing'' made of the restriction $\varphi_s$ of $s\in S$ to the set
$\{\omega\in (X_0,\mu_0)^\Gamma: \pi(s.\omega)=\pi(\omega)\}$.
Is the equivalence relation generated by $\Phi_{\pi}$  finite?
\end{question}

\section{Graphings}\label{sect:cost}

The \dfn{cost} of a p.m.p. equivalence relation $\RR$ has been introduced by Levitt \cite{Lev95}. 
It has been studied intensively in \cite{Gab98,Gab00a}. See also \cite{Kechris-Miller, Kechris=Global-aspects=2010, 2009=Furman-survey} and the popularization paper \cite{2010=Gab-what-is}.
When an equivalence relation is generated by a group action, the relations between the generators of the group introduce redundancy in the generation, and one can decrease this redundancy by using instead \dfn{partially defined isomorphisms}.

A countable family $\Phi=(\varphi_j:A_i\overset{\sim}{\to} B_j)_{j\in J}$ of measure preserving isomorphisms between Borel subsets $A_i,B_i\subset X$ is called a \dfn{graphing}. It generates a p.m.p. equivalence relation $\RR_{\Phi}$: the smallest equivalence relation such that $x\sim \varphi_j(x)$ for $j\in J$ and $x\in A_j$. Moreover, $\Phi$ furnishes a graph structure (hence the name) $\Phi[x]$ on the class of each point $x\in X$: two points $y$ and $z$ in its class are connected by an edge whenever $z=\varphi_j^{\pm 1}(y)$ for some $j\in J$. If $\RR$ is generated by a free action of $\Gamma$ and if $\Phi$ is made of isomorphisms associated with a generating set $S$ of $\Gamma$, then the graphs $\Phi[x]$ are isomorphic with the corresponding Cayley graph of $\Gamma$.
When all the graphs $\Phi[x]$ are trees, $\Phi$ is called a \dfn{treeing}.  
If it admits a generating treeing, $\RR$ is called \dfn{treeable}. See Adams \cite{Ada88, Ada90} for the first study of treed equivalence relations.

The \dfn{cost} of $\Phi$ is the number of generators weighted by the measure of their support:
$\cost (\Phi)=\sum_{ j\in J} \mu(A_j)=\sum_{j\in J} \mu(B_j)$.
The \dfn{cost of $\RR$} is the infimum over the costs of its generating graphings:
$\cost(\RR)=\inf\{\cost(\Phi): \RR=\RR_{\Phi}\}$. It is by definition an OE-invariant.
The cost of $\RR$ is $\geq 1$ when the classes are infinite \cite{Lev95}. Together with Ornstein-Weiss' theorem this gives that every p.m.p. free action of an infinite amenable group has cost $=1$. Various commutation properties in a group $\Gamma$ also entail cost $=1$ for all of its free actions. For instance when $\Gamma=G\times H$ is the product of two infinite groups and contains at least one infinite order element or $\Gamma=\SL(n,\Zmath)$, for $n\geq 3$.
It is not difficult to see that when a finite cost graphing $\Phi$ realizes the cost of $\RR_\Phi$ then $\Phi$ is a treeing.
The main results in \cite{Gab98} claim the converse: 
\thmm{If $\Phi$ is a treeing then $\cost(\RR_\Phi)=\cost (\Phi)$. 
In particular, the free actions of the free group $\FF_n$ have cost $n$}{}
In particular, free groups of different ranks cannot have OE free actions.
The cost measures the amount of information needed to construct $\RR$. It is an analogue of the \dfn{rank} of a countable group $\Gamma$, i.e. the minimal number of generators or in a somewhat pedantic formulation, the infimum of the measures $\delta(S)$ over the generating systems $S$, where $\delta$ denotes the counting measure on the group. Similarly the \dfn{cost of $\RR$} is the infimum of the measures $\nu(C)$ over the Borel subsets $C\subset \RR$ which generate $\RR$, where $\nu$ is the measure on $\RR$ introduced in section~\ref{sect: associated vN alg}, equation~(\ref{eq: def nu}) (compare Connes' Bourbaki seminar \cite{Connes-Bourbaki-04}).

In \cite{Gab00a} the notion of \dfn{free product decomposition} $\RR=\RR_1*\RR_2$ (and more generally \dfn{free product with amalgamation}  $\RR=\RR_1*_{\RR_3}\RR_2$) of an equivalence relation over subrelations is introduced  (see also \cite{Ghy95,1999=Paulin=prop-asymptot}). Of course, when $\RR$ is generated by a free action of a group, a decomposition of $\Gamma=\Gamma_1*_{\Gamma_3} \Gamma_2$ induces the analogous decomposition of $\RR=\RR_{\Gamma_1}*_{\RR_{\Gamma_3}} \RR_{\Gamma_2}$. 
The cornerstone in cost theory is the following computation:
\thmm{$\cost(\RR_1*_{\RR_3}\RR_2)=\cost(\RR_1) +\cost(\RR_2)-\cost(\RR_3)$, when $\RR_3$ is hyperfinite (possibly trivial)}{\cite{Gab00a}}
These techniques allow for the calculation of the cost of the free actions of several groups: for instance $\SL(2,\Zmath)$ ($\cost=1+1/12$), surface groups $\pi_1(\Sigma_g)$ ($\cost=2g-1$)...
In all the examples computed so far, the cost does not depend on the particular free action of the group, thus raising the following question (which proved to be related to \dfn{rank gradient} and a low-dimensional topology problem; see \cite{Abert-Nikolov-07})
(see also Question~\ref{Q:cost-beta1}):
\begin{question}[Fixed Price Problem]
Does there exist a group $\Gamma$ with two p.m.p. free actions of non equal costs?
\end{question}
Observe that both the infimum $\cost(\Gamma)$ (\cite{Gab00a}) and the supremum $\cost^*(\Gamma)$ (\cite{Abert+Weiss:2008}) among the costs of all free p.m.p. actions of $\Gamma$ are realized by some actions.

\begin{question}[Cost for Kazdhan groups]
Does there exist a Kazdhan property (T) group with a p.m.p. free action of cost $>1$?
\end{question}

In his very rich monograph \cite{Kechris=Global-aspects=2010}, Kechris studied the continuity properties of the cost function on the space of actions and proved that $\cost(\RR)>1$ for an ergodic $\RR$ forces its outer automorphism group to be Polish. 
He also introduced the topological OE-invariant $t([\RR])$, defined as the minimum number of generators of a dense subgroup of the full group $[\RR]$ and related it with the cost \cite{Kechris=Global-aspects=2010}.  When $\RR$ is generated by a free ergodic action of $\FF_n$, Miller obtained the following lower bound: $n+1\leq t([\RR])$, and \cite{Kittrell-Tsankov} proved that $t([\RR_{\hyp}])\leq 3$ and that $t([\RR])\leq 3(n+1)$.

Lyons-Pichot-Vassout \cite{Lyons-Pichot-Vassout-08} introduced the \dfn{uniform isoperimetric constant} $h(\RR)$ for p.m.p. equivalence relations, a notion similar to that for countable groups $h(\Gamma)$. They were able to obtain the purely group theoretic sharp comparison $2 \beta^{_{(2)}}_1(\Gamma)\leq h(\Gamma)$ (where $\beta^{_{(2)}}_1(\Gamma)$ is the first $\ell^2$-Betti number). 
Two complementary inequalities from \cite{Lyons-Pichot-Vassout-08,Pichot-Vassout=2009} lead to ``$2(\cost(\RR)-1)=h(\RR)$'', thus identifying two OE-invariants of apparently different nature. 
See \cite{Lyons-Peres-book} for an application of cost to percolation theory.

\section{Dimensions}\label{sect: Dimensions}

Geometric group theory studies countable groups through their actions on ``nice spaces''. Similarly, for a p.m.p. equivalence relation (it is a groupoid \cite{2000=Anantharaman-Renault=Amen-gpoids})  $\RR$ on $(X,\mu)$, one might consider its \dfn{actions} on fields of spaces $X\ni x\mapsto \Sigma_x$, or \dfn{$\RR$-field}. 
For instance, a graphing $\Phi$ defines a \dfn{measurable field of graphs} $x\mapsto \Phi[x]$, on which the natural isomorphism $\Phi[y]\simeq \Phi[z]$ for $(y,z)\in \RR_\Phi$ induces an action of $\RR_\Phi$.
The Bass-Serre theory \cite{1976=Bass=act-trees,Ser77} relates the actions of a group on trees to its free product with amalgamation decompositions (and HNN-extensions). 
Alvarez \cite{2009=Alvarez=Bass-Serre=preprint, 2009=Alvarez=Kurosh} developped an analogous theory in the framework of equivalence relations. 
For instance an equivalence relation $\RR$ acts ``properly'' on a field of trees iff $\RR$ is treeable \cite{2009=Alvarez=Bass-Serre=preprint}. He also obtained a theorem describing the structure of subrelations of a free product \cite{2009=Alvarez=Kurosh}, analogous to Kurosh's theorem. This led in \cite{Alvarez-Gaboriau-08} to the essential uniqueness of a free product decomposition $\RR=\RR_1*\cdots*\RR_n$ when the factors are \dfn{freely indecomposable} (i.e. indecomposable as a non-trivial free product) (compare \cite{IPP05,Chifan-Houdayer-08}). See also \cite{Sako-2009=ME-rigid+bi-exact} for similar results for some free products with amalgamation over amenable groups.
\begin{definition}[\cite{Alvarez-Gaboriau-08}]
A countable group is called \dfn{measurably freely indecomposable} ({\MFI}) if all its free p.m.p. actions are freely indecomposable.
\end{definition}
Examples of {\MFI} groups are provided by non-amenable groups with $\beta^{_{(2)}}_1=0$.
\begin{question}[\cite{Alvarez-Gaboriau-08}]
Produce a {\MFI} group with $\beta^{_{(2)}}_1>0$.
\end{question}

More generally, a \dfn{simplicial $\RR$-field} is a measurable field of simplicial complexes with a simplicial action of $\RR$ (see \cite{Gab02}): the space $\Sigma^{(0)}$ of $0$-cells has a Borel structure and a measurable map $\pi$ onto $X$ with countable fibers. The  cells are defined in the fibers; $\RR$ permutes the fibers; and everything is measurable. The action is \dfn{discrete} (or \dfn{smooth}, or \dfn{proper}) if it admits a measurable fundamental domain in $\Sigma^{(0)}$. For example, consider a free p.m.p. action $\Gamma\action{\alpha}(X,\mu)$ and a free action of $\Gamma$ on a (usual, countable) simplicial complex $L$. This defines a proper simplicial $\RR_{\alpha}$-action on $X\times L$ induced by the diagonal $\Gamma$-action. It is instructive to consider an OE action $\Lambda \ \action{\beta}(X,\mu)$ and to try to figure out the action of $\RR_{\beta}=\RR_{\alpha}$ on $X\times L$ once $\Gamma$ is forgotten.

The \dfn{geometric dimension}  $\geomdim (\RR)$ of $\RR$ is defined as the smallest possible dimension of a {\em proper} $\RR$-field of {\em contractible} simplicial complexes \cite{Gab02}. It is analogous to (and bounded above by) the classical geometric dimension (\cite{Bro83}) of $\Gamma$.
The \dfn{approximate dimension} \cite{Gab02} (no classical analogue) $\approxdim(\RR)$ of $\RR$ is defined 
as the smallest possible upper limit of geometric dimensions along increasing approximations of $\RR$:
$$\approxdim(\RR):=\min \{\sup (\geomdim(\RR_n))_n : (\RR_n)\nearrow \RR  \}.$$
For instance, $\geomdim (\RR)=0$ for finite equivalence relations; $\approxdim(\RR)=0$ iff $\RR$ is hyperfinite;
and $\geomdim (\RR)=1$ iff $\RR$ is treeable. Thus, quite surprisingly, surface groups admit free actions of $\geomdim=1$.
Every free action of a Kazhdan property (T) group satisfies $\approxdim= \geomdim>1$ \cite{AS90, Moo82, Gaboriau=Approx-in-prepa}.
In the following statement, $\beta^{_{(2)}}_n$ denotes the $n$-th $\ell^2$-Betti number (see section~\ref{sect: L2-Betti}).
\thmm{These dimensions satisfy:\\
-a- $\geomdim (\RR)-1\leq \approxdim(\RR)\leq \geomdim (\RR)$.\\
-b- If $\Lambda<\Gamma$ satisfies $\beta^{_{(2)}}_p(\Lambda)\not=0$, then $\geomdim(\RR_\alpha)\geq p$ for every free p.m.p. action $\Gamma\action{\alpha}(X,\mu)$. If moreover $\geomdim(\RR_\alpha)= p$, then $\beta^{_{(2)}}_p(\Gamma)\not=0$}{\cite{Gaboriau=Approx-in-prepa}}
It follows that every free action action of $\FF_{r_1}\times \cdots\times \FF_{r_p}$ ($r_j\geq 2$) (resp. $\Zmath\times \FF_{r_1}\times \cdots\times \FF_{r_p}$) has $\approxdim= \geomdim =p$ (resp. $\geomdim =p+1$).
Moreover, for every $p\geq 3$, there is a group $\Gamma_p$ with free actions $\alpha_p$ and $\beta_p$ such that $\approxdim(\RR_{\alpha_p})= \geomdim (\RR_{\alpha_p})=p$ and $\approxdim(\RR_{\beta_p})+1= \geomdim (\RR_{\beta_p})=p$.

In \cite{Dooley-Golodets-2009}, Dooley-Golodets study the behavior of the dimension $\geomdim$ under finite extensions.
The notion of \dfn{measurable cohomological dimension} introduced in \cite{Sauer-Thom=spectral-seq=L2=2007} has some similarity with the geometric dimension.

\section{$L^2$-Betti numbers}\label{sect: L2-Betti}

The $\ell^2$-Betti numbers of cocompact group actions on manifolds were introduced by Atiyah \cite{Ati76} in terms of the heat kernel. Connes \cite{Connes-79} defined them for measured foliations. Cheeger-Gromov \cite{CG86} introduced $\ell^2$-Betti numbers $\beta^{_{(2)}}_n(\Gamma)\in [0,\infty]$, $n\in \Nmath$, for arbitrary countable groups $\Gamma$. In \cite{Gab02} the $L^2$-Betti numbers $\beta^{_{(2)}}_n(\RR)\in [0,\infty]$, $n\in \Nmath$ are defined for p.m.p. equivalence relations $\RR$, by using proper simplicial $\RR$-fields (see section~\ref{sect: Dimensions}). In any case, the definitions rely on the notion of generalized von Neumann dimension, expressed as the trace of certain projections. 
One of the main results in \cite{Gab00b,Gab02} is the invariance of the $\beta^{_{(2)}}_n(\Gamma)$ under orbit equivalence.
\thmm{If $\RR_{\Gamma}$ is generated by a free p.m.p. action of $\Gamma$, then $\beta^{_{(2)}}_n(\RR_\Gamma)=\beta^{_{(2)}}_n(\Gamma)$ for every $n\in \Nmath$}{\cite{Gab02} }

The inequality $\cost(\Gamma)\geq\beta^{_{(2)}}_1(\Gamma)-\beta^{_{(2)}}_0(\Gamma)+1$ proved in \cite{Gab02} is an equality in all cases where the computations have been achieved, thus leading to the question: 
\begin{question}[Cost vs first $\ell^2$-Betti number]\label{Q:cost-beta1}
Is there an infinite countable group with  $\cost(\Gamma)>\beta^{_{(2)}}_1(\Gamma)+1$ ?
\end{question}
The following compression formula was a key point in various places notably when studying ``self-similarities'' (the ``fundamental group'', see \cite{Pop06}) and measure equivalence (see section~\ref{sect:Measure equivalence}).
\thmm{The $L^2$-Betti numbers of $\RR$ and of its \dfn{restriction} to a Borel subset $A\subset X$ meeting all the classes satisfy: $\beta^{_{(2)}}_n(\RR)=\mu(A)\beta^{_{(2)}}_n(\RR\vert A)$}{\cite{Gab02}}
It follows that lattices in a common locally compact second countable group have proportional $\ell^2$-Betti numbers.

In \cite{BG04}, $L^2$-Betti numbers for profinite actions are used to extend L\"uck's approximation theorem \cite{Luc94b} to non-normal subgroups.
We refer to the book \cite{Luc02} for information about $\ell^2$-Betti numbers of groups and for an alternative approach to von Neumann dimension.
See \cite{Sauer=Betti-groupoids=2005,Sauer-Thom=spectral-seq=L2=2007,Thom=L2-inv+rank-metric=2008} for extension of  $\beta^{_{(2)}}_n(\RR)$ to measured groupoids, and several computations using L\"uck's approach (\cite{Neshveyev-Rustad=def-L2=2009} proves that the various definitions coincide).

Very interesting combinatorial analogues of the cost and $\beta^{_{(2)}}_1$ have been introduced by Elek \cite{Elek=combinatorial-cost=2007} in a context of sequences of finite graphs.

\section{Measure equivalence}\label{sect:Measure equivalence}

Two groups $\Gamma_1$ and $\Gamma_2$ are \dfn{virtually isomorphic} if there exist $F_i \triangleleft \Lambda_i<\Gamma_i$ 
such that $\Lambda_1/F_1\simeq \Lambda_2/F_2$, where $F_i$ are finite groups, and $\Lambda_i$ has finite index in $\Gamma_i$.
This condition is equivalent with: $\Gamma,\Lambda$ admit commuting actions on a set $\Omega$ such that each of the actions $\Gamma\acting\Omega$ and $\Lambda\acting\Omega$ has finite quotient set and finite stabilizers.

A finite set admits two natural generalizations, a topological one (compact set) leading to \dfn{geometric group theory} and a measure theoretic one (finite measure set) leading to \dfn{measured group theory}. 
\begin{definition}[\cite{Gro93}]
Two countable groups $\Gamma_1$ and $\Gamma_2$ are \dfn{measure equivalent (ME)} (\textbf{notation:} $\Gamma_1\ME \Gamma_2$) if there exist commuting actions of $\Gamma_1$ and $\Gamma_2$, that are (each) measure preserving, free, and with a finite measure fundamental domain, on some standard (infinite) measure space $(\Omega, m)$. 
\end{definition}
The ratio $[\Gamma_1:\Gamma_2]_\Omega:=m(\Omega/\Gamma_2)/m(\Omega/\Gamma_1)$ of the measures of the fundamental domains is called the \dfn{index} of the \dfn{coupling} $\Omega$.
The typical examples, besides virtually isomorphic groups, are lattices in a common (locally compact second countable) group $G$ with its Haar measure, acting by left and right multiplication. 

The topological analogue was shown to be equivalent with \dfn{quasi-isometry (QI)} between finitely generated groups \cite{Gro93}, thus raising \dfn{measured group theory} (i.e. the study of groups up to ME) to parallel \dfn{geometric group theory}. 
See \cite{Fur99a} for the basis in ME and the surveys \cite{Gab05a,Shalom=survey=2005,2009=Furman-survey} for more recent developments.
Measure equivalence and orbit equivalence are intimately connected by considering the relation between the quotient actions $\Gamma_1\acting \Omega/\Gamma_2$ and $\Gamma_2\acting \Omega/\Gamma_1$. In fact two groups are ME iff they admit SOE free actions.

\begin{definition}[Stable Orbit Equivalence]\label{def: SOE}
Two p.m.p. actions of $\Gamma_i\acting(X_i,\mu_i)$ are {\bf stably orbit equivalent (SOE)}
if there are Borel subsets $Y_i\subset X_i$,  ${i=1,2}$ which meet almost every orbit of $\Gamma_i$ and
a measure-scaling isomorphism $f:Y_1\to Y_2$ s.t. 
\begin{center}$f(\Gamma_1.x\cap Y_1)=\Gamma_2.f(x)\cap Y_2$ \ \ a.e.\end{center}
The \dfn{index}  or \dfn{compression constant} of this SOE $f$ is $[\Gamma_1:\Gamma_2]_f=\frac{\mu(Y_2)}{\mu(Y_1)}$.
\end{definition}

The state of the art ranges from quite well understood ME-classes to mysterious and very rich examples.
For instance, the finite groups obviously form a single ME-class. The infinite amenable groups form a single ME-class \cite{OW80}.
The ME-class of a lattice in a center-free simple Lie group $G$ with real rank $\geq 2$ (like $\SL(n,\Rmath)$, $n\geq 3$) consists in those groups that are virtually isomorphic with a lattice in $G$ \cite{Fur99a}.
If $\Gamma$ is a non-exceptional mapping class group, its ME-class consists only in its virtual isomorphism class \cite{Kida-08}.
Kida extended this kind of result to some amalgamated free products (see \cite{Kida-2009}).

On the opposite, the ME-class of the (mutually virtually isomorphic) free groups $\FF_r$ ($2\leq r <\infty$) contains the free products $*_{i=1}^{r} A_i$ of infinite amenable groups, surface groups $\pi_1(\Sigma_g)$ ($g\geq 2$), certain branched surface groups \cite{Gab05a}, elementarily free groups \cite{Bridson-Tweedale-Wilton-2007}... and is far from being understood. 
Being ME with a free group is equivalent to admitting a free p.m.p. treeable action \cite{Hjo06}.

There is a considerable list of ME-invariants (see \cite{Gab05a} and the references therein). For instance Kazhdan property (T), Haagerup property, the \dfn{ergodic dimension} (resp. \dfn{approximate ergodic dimension}) defined as the infimum of the geometric (resp. approximate) dimension among all the free p.m.p. actions of $\Gamma$, the sign of the Euler characteristic (when defined), the Cowling-Haagerup invariant, belonging to the classes $\mathcal{C}_{\mathrm{reg}}, \mathcal{C}$.
Recently exactness (see \cite{Brown-Ozawa08}) and belonging to the class $\mathcal{S}$ of Ozawa \cite{Sako-2009=class-S-ME} were proved to be ME-invariants. There are also numerical invariants which are preserved under ME modulo multiplication by the index: $\cost(\Gamma)-1$, the $\ell^2$-Betti numbers $(\beta^{_{(2)}}_n(\Gamma))_{n\in \Nmath}$ \cite{Gab02}.

ME is stable under some basic constructions: \\
(a) if $\Gamma_i\ME\Lambda_i$ for $i=1,\cdots,n$ then $\Gamma_1\times\cdots\times\Gamma_n\ME \Lambda_1\times\cdots\times\Lambda_n$\\
(b) if $\Gamma_i\ME\Lambda_i$ with index 1, then $\Gamma_1*\cdots*\Gamma_n\ME \Lambda_1*\cdots*\Lambda_n$ (with index 1).

Some papers study when the converse holds \cite{MS06,IPP05,Chifan-Houdayer-08,Alvarez-Gaboriau-08}. One has of course to impose some irreducibility conditions on the building blocks, and these conditions have to be strong enough to resist the measurable treatment. 
These requirements are achieved \\
(a) (for direct products) if the $\Gamma_i, \Lambda_i$ belong to the class $\mathcal{C}_{\mathrm{reg}}$ of \cite{MS06} (for instance if they are non-amenable non-trivial free products): the non-triviality of the bounded cohomology ${{\rm H}^2_{\rm b}}(\Gamma,\ell^2(\Gamma))$ is an ME-invariant preventing $\Gamma$ to decompose (non-trivially) as a direct product;\\
(b) (for free products) if the $\Gamma_i, \Lambda_i$ are \MFI (for instance if they have $\beta^{_{(2)}}_1=0$ and are non-amenable) \cite{Alvarez-Gaboriau-08}: they are not ME with a (non-trivial) free product. We prove for instance:
\begin{theorem}[\cite{Alvarez-Gaboriau-08}]
If $\Gamma_1*\cdots*\Gamma_n\ME \Lambda_1*\cdots*\Lambda_p$, where both the $\Gamma_i$'s and the $\Lambda_j$'s belong to distinct ME-classes and are \MFI, then $n=p$ and up to a permutation of the indices $\Gamma_i\ME\Lambda_i$.
\end{theorem}
See also \cite{IPP05,Chifan-Houdayer-08} when the groups have Kazhdan property (T), or are direct products, under extra ergodicity hypothesis. The delicate point of removing ergodicity assumptions in \cite{Alvarez-Gaboriau-08} was achieved by using \cite{2009=Alvarez=Kurosh}.

Similar ``deconstruction'' results were obtained by Sako \cite{Sako-2009=ME-rigid+bi-exact} for building blocks made of direct products of non-amenable exact groups when considering free products with amalgamation over amenable subgroups or by taking wreath product with amenable base.

Refinements of the notion of ME were introduced in \cite{Shalom04,Thom=bound-coh-L2=2009,Lueck-Sauer-Wegner=preprint-L1-ME=2009} or by Sauer and Bader-Furman-Sauer.

\section{Non-orbit equivalent actions for a given group}\label{sect:How many actions}
In this section, we only consider ergodic free p.m.p. actions $\Gamma\action{\alpha}(X,\mu)$ of infinite countable groups and the associated orbit equivalence relations $\RR_\alpha$.
Ornstein-Weiss' theorem \cite{OW80} implies that amenable groups all produce the same relation, namely $\RR_\hyp$. What about non-amenable groups? How many non-OE actions for a given group? Most of the OE-invariants depend on the group rather than on the action, and thus cannot distinguish between various actions of the group.
However, for non-Kazhdan property (T) groups, Connes-Weiss \cite{Connes-Weiss=1980} produced two non-OE actions distinguished by strong ergodicity (see section~\ref{sect:strong ergodicity}). 
And along the years, various rigidity results entailed some specific families of groups to admit continously\footnote{This is an upper bound since $\mathrm{Card}(\Aut([0,1],\Leb))=2^{\aleph_0}$.} many non-OE actions (see for instance \cite{1981=Bezuglyi-Golodets,Zimmer-Erg-th-ss-gp-Book-84, GG88, MS06,Popa-1-cohomology-2006,Popa-OE-superrig-2007}).

We briefly describe below the crucial steps on the route toward the general solution:
\thmm{Every non-amenable group admits continuously many orbit inequivalent free ergodic p.m.p. actions}{\cite{Ioa07,Epstein-07}}
The first step was made by Hjorth \cite{Hjorth-conv-Dye-05} when, within the circle of ideas from Connes \cite{Co80} and Popa \cite{Popa=correspondences}, he obtained the result for Kazhdan property~(T) groups. Roughly speaking, a pair of OE actions $\alpha$ and $\beta$ defining the same equivalence relation $\RR$ gives a diagonal action ($\gamma.(x,y)=(\gamma._\alpha x, \gamma._\beta y)$) on $\RR$  and thus a unitary representation on $L^2(\RR,\nu)$. When considering uncountably many OE actions, a separability argument shows that the characteristic function $\mathbf{1}_D$ of the diagonal is sufficiently almost invariant for some pair of actions. Now, an invariant vector near $\mathbf{1}_D$, which is given by property (T), delivers a conjugacy between the actions. There exists a continuum of pairwise non-conjugate actions, and by the above the OE-classes in this continuum are countable.

The next step was the analogous theorem for the prototypical non-property~(T), non-rigid group, namely the free groups and some free products \cite{GP05}. It lay again within the same circle of ideas but there, rigidity was obtained through Popa's \dfn{property (T) relative to the space} (see section~\ref{sec:Rel T}). 

Then Ioana \cite{Ioa07} extended it to all groups containing a copy of $\FF_2$. For this, he introduced a weak version of property (T) relative to the space and used a general construction called \dfn{co-induction}\footnote{Co-induction is the classical right adjoint of restriction. Its measure theoretic version was brought to my attention by Sauer and used in \cite{Gab05a}, but it probably first appeared in preliminary versions of \cite{DGRS=2008}.}.

Eventually, Epstein obtained the theorem in full generality \cite{Epstein-07}. For this she had to generalize the co-induction construction to the setting provided by Gaboriau-Lyons' measurable solution to \dfn{von Neumann's problem} (see below). Moreover, Ioana extended Epstein's result from 
orbit inequivalent to von Neumann inequivalent actions \cite{Ioa07}.

When von Neumann introduced the notion of amenability \cite{vN29}, he observed that a countable group containing a copy of $\FF_2$ cannot be amenable. The question of knowing whether every non-amenable countable group has to contain a copy of $\FF_2$, known as \dfn{von Neumann's problem}, was answered in the negative by Ol$'${\v{s}}anski{\u\i} \cite{OlShanskii=1980}. In the measurable framework, offering much more flexibility, the answer is somewhat different: 
\thmm{For any non-amenable countable group $\Gamma$, the orbit equivalence relation of the Bernoulli shift action $\Gamma\acting ([0,1],\Leb)^{\Gamma}$ contains a subrelation generated by a free ergodic p.m.p. action of $\FF_2$}{\cite{Gaboriau-Lyons-2009}}
In the terminology of \cite{Monod=ICM-2006}, there is a randembedding of $\FF_2$ in any non-amenable group.
The proof uses percolation theory on graphs and \cite{Haggstrom-Peres=1999,Lyons-Schramm=1999,PS00, Gab05,Hjo06}. The following general question remains open:

\begin{question}
Does every ergodic non-hyperfinite p.m.p. equivalence relation contain a (treeable) subrelation of cost $>1$?
\end{question}

\section{Relative property (T)}\label{sec:Rel T}
In his seminal paper \cite{Kaz67} on property (T), Kazhdan implicitly\footnote{This was made explicit in \cite{1983=Margulis=fin-add-meas-Eucli-sp}.} introduced the notion of property (T) relative to a subgroup $\Lambda<\Gamma$. In particular, a group always has property (T) relative to its ``unit subgroup'' $\{1\}<\Gamma$. When considering a groupoid like $\RR$, its space of units $(X,\mu)$ (and its ``relative representation theory'') is much more complicated. 
The introduction by Popa \cite{Pop06} of the fruitful notion of \dfn{property (T) relative to the space $(X,\mu)$} (also simply called \dfn{rigidity}) allowed him to solve some long standing problems in von Neumann algebras. In fact, the definition involves a pair of von Neumann algebras $B\subset M$ (for instance $L^\infty(X,\mu)\subset L(\RR)$) and parallels the analogous notion for groups, in the spirit of Connes-Jones \cite{CJ85}.

The typical example is provided by the standard action of $\SL(2,\Zmath)$ and its non-amenable subgroups $\Gamma$ (for instance free groups $\FF_r$, $r\geq 2$) on $\Tmath^2$.
Notice Ioana's result that in fact every ergodic non-amenable subrelation of $\RR_{\SL(2,\Zmath)\acting \Tmath^2}$ still has property (T) relative to the space $\Tmath^2$ \cite{ioana-2009}.
The property (T) relative to the space $(X,\mu)$ comes from the group property (T) of $\Zmath^2\rtimes \Gamma$ relative to the subgroup $\Zmath^2$, via viewing $\Zmath^2$ as the Pontryagin dual of $\Tmath^2$. 
This property (never satisfied by standard Bernoulli shifts) entails several rigidity phenomena (see for instance \cite{Pop06,IPP05,GP05}). More examples come from \cite{Val05,Fer06} and they all involve some arithmeticity. This led Popa to ask for the class of groups admitting such a free p.m.p. action with property (T) relative to the space. T\"ornquist \cite{Tor06} ensures that the class is stable under taking a free product with any countable group.  More generally, \cite{Gab08} shows that the class contains all the non-trivial free products of groups $\Gamma=\Gamma_1*\Gamma_2$: in fact $\RR_{\Gamma_1}$ and $\RR_{\Gamma_2}$ may be chosen to be conjugate with any prescribed free $\Gamma_i$-action and the arithmeticity alluded to is hidden in the way they are put in free product. This leads, using ideas from \cite{PV08} to (plenty of) examples of $\RR_\Gamma$ with trivial outer automorphism group, in particular the first examples for free $\FF_2$-actions \cite{Gab08}.
Ioana \cite{Ioa07} proved that every non-amenable group admits a free p.m.p. action satisfying a weak form of the above property, enough for various purposes, see section~\ref{sect:How many actions}.

\section{Some rigidity results}\label{sec:some rigidity results}
We have three notions of equivalence between free p.m.p. actions:
\begin{equation*}
(\Gamma_1\action{\alpha_1}X_1 \Conj \Gamma_2\action{\alpha_2}X_2)\ \Longrightarrow\ (\Gamma_1\action{\alpha_1}X_1 \OrbEqu \Gamma_2\action{\alpha_2}X_2)\ \Longrightarrow\ (\RR_{\alpha_1} \vNEqu \RR_{\alpha_2}).
\end{equation*}
\dfn{Rigidity} phenomena consist ideally in situations where (for free actions) some implication can be reversed, or more generally when a big piece of information of a stronger nature can be transferred through a weaker equivalence.
Zimmer's pioneering work (see \cite{Zimmer-Erg-th-ss-gp-Book-84}) inaugurated a series of impressive results of rigidity for the first arrow 
($\Conj {\dashleftarrow} \OrbEqu$), made possible by the introduction in OE theory and in operator algebras of new techniques borrowed from diverse mathematical domains, like algebraic groups, geometry, geometric group theory, representation theory or operator algebras. These rigidity results for $\Gamma_1\action{\alpha_1}X_1$ take various qualifications according to whether an OE hypothesis entails\\
 -- \dfn{strong OE rigidity}: conjugacy under some additionnal hypothesis about the mysterious action $\Gamma_2\action{\alpha_2}X_2$, or even  \\
-- \dfn{OE superrigidity}: conjugacy of the actions with no hypothesis at all on the target action.
\\
These notions are \dfn{virtual} when they happen only up to finite groups (see \cite{Fur99b} for precise definitions).

To give some ideas we simply evoke a sample of some typical and strong statements far from exhaustiveness or full generality. 
\thmm{Any free action that is OE with the standard action $\SL(n,\Zmath)\acting\Tmath^n$ for $n\geq 3$, is virtually conjugate with it}{\cite{Fur99b}}
This is more generally true for lattices in a connected, center-free, simple, Lie group of higher rank, and for ``generic'' actions (see \cite{Fur99b}).
Monod-Shalom \cite{MS06} obtained strong OE rigidity results when $\Gamma_1$ is a direct product of groups in $\mathcal{C}_{\mathrm{reg}}$, under appropriate ergodicity assumptions on both sides.
See also Hjorth-Kechris \cite{Hjorth-Kechris=Rigid+products-05} for rigidity results about actions of products, where the focus is more on Borel reducibility. Kida's results \cite{Kida-2008-OE} consider actions of mapping class groups of orientable surfaces and their direct products. He also obtains very strong rigidity results for certain amalgamated free products \cite{Kida-2009}.
A series of ground breaking results in von Neumann algebras obtained by Popa 
 \cite{Pop06, Popa-strong-rigid-I-2006, Popa-strong-rigid-II-2006, Popa-OE-superrig-2007, Popa-superrigid-spect-gap-2008} 
 and his collaborators \cite{PS07, IPP05, Popa-Vaes-2008, Popa-Vaes=superrig-SLn=2008,PV08, Ioana-08, Popa-Vaes=fund-gp=2008-preprint, Popa-Vaes=W-rigidity=2009-preprint} (see \cite{Vaes=Bourbaki=2007} for a review) dramatically modified the landscape. On the OE side, these culminated in Popa's cocycle superrigidity theorems, that imply several impressive OE superrigidity corollaries, for instance:
\thmm{Assume that $\Gamma$ is either an infinite ICC Kazhdan property (T) group or is the product of two infinite groups $H\times H'$ and has no finite normal subgroup.
Then any free action that is orbit equivalent with the Bernoulli shift $\Gamma\acting (X_0,\mu_0)^{\Gamma}$ is conjugate with it}{\cite{Popa-OE-superrig-2007, Popa-superrigid-spect-gap-2008}}
See Furman's ergodic theoretical treatment and generalizations \cite{Furman-08} for the Kazhdan property (T) case. 
In the opposite direction, Bowen obtained some surprising non-rigidity results \cite{bowen-2009-coinduced,bowen-2009-Bernoul-free-gp} showing for instance that all the Bernoulli shifts of the free groups $\FF_r$, $2\leq r<\infty$ are mutually SOE (see Def.~\ref{def: SOE}).

\medskip
As it follows from  \cite{Singer-1955,FM77b}, being able to reverse the second arrow ($\OrbEqu {\dashleftarrow} \vNEqu$) essentially amounts to being able to uniquely identify the Cartan subalgebra inside $L(\RR)$, i.e. given two Cartan subalgebras $A_1, A_2$ in $L(\RR_1)\simeq L(\RR_2)$, being able to relate them through the isomorphism. Such results are qualified \dfn{vNE rigidity} or \dfn{$W^*$-rigidity}.
 The starting point is Popa's breakthrough \cite{Pop06} where a uniqueness result is obtained under some hypothesis on both $A_1$ and $A_2$ (and this was enough to solve long standing problems in von Neumann algebras). See also \cite{IPP05, Chifan-Houdayer-08} for this kind of strong statements under various quite general conditions. We refer to the surveys \cite{Popa=ICM-2007,Vaes=Bourbaki=2007,Vaes=ICM=2010} for the recent developments in vNE or $W^*$-rigidity.
 However, after a series of progresses (see for instance \cite{Ozawa-Popa-08, ozawa-Popa-II=2008-b,Ioana-08,Peterson=L2-rigid=2009, Popa-Vaes=W-rigidity=2009-preprint,Peterson-2010=preprint}), the most recent achievement is:
 \thmm{If a free action of a group is von Neumann equivalent with the standard Bernoulli shift action of an ICC Kazhdan property (T) group, then the actions are in fact conjugate}{\cite{Ioana=W-rigidity=2010=preprint}}

\section{Some further OE-invariants}

In order to distinguish treeable Borel equivalence relations, Hjorth introduced a technique preventing a p.m.p. equivalence relation from being OE with a profinite one \cite{Hjo06}. Then Kechris and Epstein-Tsankov isolated representation-theoretic properties (i.e. in terms of the Koopman representation) leading to strong forms of non-profiniteness; see \cite{Kechris-2005-modular,Epstein-Tsankov=modular=2010}.

Elek-Lippner introduced the \dfn{sofic} property for equivalence relations. It is satisfied by profinite actions, treeable equivalence relations and Bernoulli shifts of sofic groups \cite{Elek-Lippner=sofic-rel=2010}. They also proved that the associated von Neumann algebra satisfies the Connes' embedding conjecture.

\medskip
\textbf{Acknowledgements} I'm grateful to A. Alvarez, C. Houdayer and J. Melleray for their comments.

\bibliographystyle{alpha}
\def\cprime{$'$} \def\cprime{$'$} \def\cprime{$'$}

\bigskip
{\footnotesize \hskip-\parindent Damien Gaboriau\\
Unit\'e de Math\'ematiques Pures et Appliqu\'ees\\
Universit\'e de Lyon, CNRS, Ens-Lyon\\
46, all\'ee d'Italie \\
69364 Lyon cedex 7, France\\
{gaboriau@umpa.ens-lyon.fr}}

\end{document}